\documentclass[11pt]{article}

\usepackage[utf8]{inputenc} 
\usepackage[T1]{fontenc}    
\usepackage{url}            
\usepackage{booktabs}       
\usepackage{amsfonts}       
\usepackage{nicefrac}       
\usepackage{microtype}      
\usepackage{lipsum}
\usepackage{amsmath}
\usepackage{amsthm}
\usepackage{graphicx}
\usepackage{subcaption}
\usepackage{amssymb}

\newcommand{\kk}{\kappa}
\newcommand{\bb}{\beta}

\newcommand{\C}{\mathbb{C}}

\newcommand{\N}{\mathbb{N}}
\newcommand{\Z}{\mathbb{Z}}

\newcommand\mydots{\text{\makebox[1em][c]{.\hfil.\hfil.}}}
\providecommand{\keywords}[1]
{
  \small	
  \textbf{Keywords:} #1
}

\title{Zeros of Bessel cross-products coming from oblique derivative boundary value problems}

\author{
  Stanislav~Budzinskiy \\
  Faculty of Computational Mathematics and Cybernetics\\
  Lomonosov Moscow State University
}
\date{}

\begin{document}
\maketitle

\begin{abstract}
The paper is devoted to (combinations of) Bessel cross-products that arise from oblique derivative boundary value problems for the Laplacian in a circular annulus. We show that like their Neumann-Laplacian counterpart (and unlike the Dirichlet-Laplacian), they possess two kinds of zeros: those that can be derived by McMahon series and diverge to infinity in the limit, and exceptional ones that remain finite. For both cases we find asymptotic expressions for a fixed oblique angle and vanishing thickness of the annulus. We further present plots of numerically computed zeros and discuss their behaviour when the oblique angle changes and the thickness remains fixed.
\end{abstract}


\section{Introduction}
Let us consider Bessel functions of the first and second kinds---$J_\nu$ and $Y_\nu$, respectively---of a complex argument $z \in \C$ and real nonnegative order $\nu\geqslant 0$. These functions appear naturally in PDEs in circular domains and their zeros define the eigenvalues of the Laplace operator. For example, zeros $\lbrace j_{\nu,n} \rbrace_{n\in\N}$ of $J_\nu$ come from the Dirichlet-Laplacian in a disc and $\lbrace j'_{\nu,n}\rbrace_{n\in\N}$ of $J'_\nu$ are connected with the corresponding Neumann-Laplacian. This explains why properties of zeros of Bessel functions (of both kinds) and their derivatives were extensively studied \cite{LaforgiaNataliniZeros2007, KerimovStudies2014}.  

When the domain is no longer a disc but an annulus, there are two boundary conditions to meet and so we are interested in the so-called cross-products of Bessel functions \cite{CochranRemarks1964}. Leaving mixed problems aside, we mention
\begin{equation}\label{cp_d}
	f_{\nu}(\kk, z) = J_{\nu}(z)Y_{\nu}(\kk z) - J_{\nu}(\kk z)Y_{\nu}(z) 
\end{equation}
\noindent and 
\begin{equation}\label{cp_n}
	g_{\nu}(\kk, z) = J'_{\nu}(z)Y'_{\nu}(\kk z) - J'_{\nu}(\kk z)Y'_{\nu}(z) 
\end{equation}
that are related to the Dirichlet and Neumann boundary conditions, respectively. Here, $\kk > 0$ is the thickness parameter. The zeros of these functions received plenty of attention with a view towards waveguides \cite{WaidronTheory1957, CochranPecinaMode1966}; they are also used to compute the Pleijel constant for planar annuli \cite{BobkovAsymptotic2019}. Note that we are being ambiguous saying "zeros": we mean $z$-zeros while there are also $\nu$-zeros that are of interest \cite{CochranASYMPTOTIC1966, CochranAnalyticity1966}.

The $z$-zeros of both \eqref{cp_d} and \eqref{cp_n} are well-studied. These functions are even, and it is known that each has a countable set of simple real $z$-zeros. When $\kk \to 1$, the positive ones can be obtained as McMahon series \cite{CochranRemarks1964,McMahonRoots1894}
\begin{equation}\label{mcmahon}
z_s = \frac{s\pi}{\kappa - 1} + p\left( \frac{s\pi}{\kappa - 1} \right)^{-1} + (q - p^2) \left( \frac{s\pi}{\kappa - 1} \right)^{-2} + \ldots, \quad s \in \N.
\end{equation}
For the zeros $z_{s}^{\rm D}$ of \eqref{cp_d} we have
\begin{equation*}
p = \frac{4\nu^2 - 1}{8\kappa}, \quad q = \frac{(\kappa^{3} - 1) (16\nu^4 - 104\nu^2 + 25)}{384 \kappa^{3} \left(\kappa - 1\right)},
\end{equation*}
while those of \eqref{cp_n}, $z_{s}^{\rm N}$, are defined by
\begin{equation*}
p = \frac{4\nu^2 + 3}{8\kappa}, \quad q = \frac{(\kappa^{3} - 1) (16\nu^4 + 184\nu^2 - 63)}{384 \kappa^{3} \left(\kappa - 1\right)}.
\end{equation*}
But numerical computations showed \cite{TruellConcerning1943} that there was one more zero of \eqref{cp_n} that was not captured by \eqref{mcmahon} . Later its series representation was discovered \cite{BuchholzBesondere1949} 
\begin{equation*}
\frac{\nu}{z_{0}^{\rm N}} = \sqrt{\kappa}\left[ 1 + \frac{(\kappa-1)^2}{12\kappa} + \frac{(8\nu^2-3)(\kappa-1)^4}{480\kappa^2} + \ldots\right], \quad \nu \neq 0.
\end{equation*}
And another form is known \cite{GottliebEigenvalues1985}
\begin{equation*}
z_{0}^{\rm N} = \nu\left[ 1 - \frac{1}{2}(\kk - 1) + \frac{7}{24}(\kk - 1)^2 + \ldots \right], \quad \nu \neq 0.
\end{equation*}
This zero is "exceptional" in that it stays finite in the limit.

In this paper we consider Bessel cross-products that stem from an oblique derivative boundary value problem in an annulus with a constant oblique angle:
\begin{multline*}
	g_{\nu}(\bb,\kk, z) = [J_\nu'(z)Y_\nu'(\kappa z) - J_\nu'(\kappa z)Y_\nu'(z)] - i\frac{\bb \nu}{z}[J_\nu(z)Y_\nu'(\kappa z) - J_\nu'(\kappa z)Y_\nu(z)] -\\- i\frac{\bb \nu}{\kappa z}[J_\nu'(z)Y_\nu(\kappa z) - J_\nu(\kappa z)Y_\nu'(z)] - \frac{\bb^2\nu^2}{\kappa z^2}[J_\nu(z)Y_\nu(\kappa z) - J_\nu(\kappa z)Y_\nu(z)] = 0.
\end{multline*}
Here, $\bb$ stands for the tangent of the oblique angle. To our knowledge, similar functions were only dealt with in \cite{Martineknote1968}, where $\nu$-zeros were studied for a combination of cross-products with real coefficients arising from elasticity. To simplify notation we introduce
\begin{equation*}
	G^{m,k}_{\nu}(\kk, z) = J_\nu^{(m)}(z)Y_\nu^{(k)}(\kappa z) - J_\nu^{(k)}(\kappa z)Y_\nu^{(m)}(z).
\end{equation*}
Since $g_{\nu}$ satisfies
$$g_{\nu}(\bb,1/\kk, \kk z) = -g_{\nu}(\bb,\kk, z),\quad g_{\nu}(\bb,\kk, -z) = g_{\nu}(\bb,\kk, z),\quad \overline{g_{\nu}(\bb,\kk, z)} = g_{\nu}(-\bb,\kk, \overline{z}),$$
we can restrict ourselves to searching for $z$-zeros of
\begin{equation}\label{g}
	g_{\nu}(\bb,\kk, z) = G^{1,1}_{\nu}(\kk, z) - i\frac{\bb \nu}{z} G^{0,1}_{\nu}(\kk, z) - i\frac{\bb \nu}{\kappa z} G^{1,0}_{\nu}(\kk, z) - \frac{\bb^2\nu^2}{\kappa z^2} G^{0,0}_{\nu}(\kk, z)
\end{equation}
with $\kk > 1$, $\bb \geqslant 0$, and $\mathrm{Re}z \geqslant 0$.

\section{Exceptional zero}
In this section we will follow the approach of \cite{Gottliebexceptional1985,GrebenkovAnalytical2008}.
As we are going to deal with a finite-in-the-limit zero, we represent it as a regular perturbation series in $\varepsilon = \kk - 1$
\begin{equation*}
	z = z_0 + z_1 \varepsilon + z_2 \varepsilon^2 + \ldots
\end{equation*}
Now we need to expand $C^{(p)}(z)$ and $C^{(p)}(\kk z)$ in powers of $\varepsilon$, where $C(z)$ is analytic at $z_0$. Since powers of power series can be expressed in terms of ordinary partial Bell polynomials
\begin{equation*}
	(z - z_0)^{n} = \sum_{k = n}^{\infty} \hat{B}_{k,n}(z_1,\mydots,z_{k-n+1})\varepsilon^k,
\end{equation*}
we find 
\begin{equation*}
	C^{(p)}(z) = \sum_{k = 0}^{\infty} \left( \sum_{j = 0}^{k} C^{(p + j)}(z_0) u_{k,j} \right) \varepsilon^k, \quad
    u_{k,j} = \frac{1}{j!}\hat{B}_{k,j}(z_1,\mydots,z_{k-j+1}).
\end{equation*}
The series for $C^{(p)}(\kk z)$ is slightly more complicated
\begin{equation*}
	C^{(p)}(\kk z) = \sum_{m = 0}^{\infty} \left( \sum_{n = 0}^{m} C^{(p + n)}(z_0) v_{m,n} \right) \varepsilon^m
\end{equation*}
with
\begin{equation*}
    v_{m,n} = \frac{1}{n!}\sum_{k = 0}^{n} \binom{n}{k} \sum_{t = 0}^{m - n} \hat{B}_{t+k,k}(z_1,\mydots,z_{t+1}) \hat{B}_{m-t-k,n-k}(z_0,\mydots,z_{m-n-t})
\end{equation*}
This leads to
\begin{equation*}
	G^{p,q}_{\nu}(\kk, z) = \sum_{l = 0}^{\infty} \left( \sum_{r = 0}^{l} \sum_{j = 0}^{r} \sum_{n = 0}^{l - r} u_{r,j} v_{l-r,n} G^{p+j,q+n}_{\nu}(1, z_0) \right) \varepsilon^{l}.
\end{equation*}
Every $G^{p,q}_{\nu}(1, z_0)$ can be expressed in terms of $G^{0,1}_{\nu}(1, z_0) = 2 / (\pi z_0)$ (which is the Wronskian \cite[Eq. 9.1.16]{AbramowitzStegunHandbook2013}) since $G^{p,q}_{\nu}(\kk, z_0) + G^{q,p}_{\nu}(\kk, z_0) = 0$. The resulting series thus takes the form
\begin{equation*}
	\kk z^2 g_{\nu}(\bb,\kk,z) = \frac{2}{\pi z_0} \left( a_0 + a_1(z_0)\varepsilon + a_2(z_0, z_1)\varepsilon^2 + a_3(z_0, z_1, z_2)\varepsilon^3 + \ldots \right),
\end{equation*}
where
\begin{align*}
	a_0 &= 0,\\
	a_1 &= z_{0} \left( z_{0}^{2} - \nu^{2}(\beta^2 + 1) \right),\\
	a_2 &= \frac{z_{0}}{2} \left(\nu^2 (\beta^{2} + 1) + z_{0}^{2} + 4 z_{0} z_{1}\right),\\
	a_3 &= - \frac{\nu^2 z_0(\beta^{2} + 1)}{6} \left(z_{0}^{2} -\nu^{2} - 2\right) + \frac{i \beta}{3} \nu z_{0}^{3} + \frac{\nu^{2} z_{0}^{3}}{6} + 2 \nu^{2} z_{1} - \frac{z_{0}^{5}}{6} - z_{0}^{2} z_{1} + 2 z_{0}^{2} z_{2} + z_{0} z_{1}^{2}.
\end{align*}
Since powers of $\varepsilon$ are independent, it is left to solve equations $a_k = 0$ one by one whence we find
\begin{equation}\label{exceptional}
	z = \nu\sqrt{\bb^2 + 1} \left[ 1 - \frac{\varepsilon}{2} - \left( \frac{5\bb^2 - 7}{\bb^2 + 1} + 4 i \bb \nu \right)\frac{\varepsilon^2}{24} + \ldots \right].
\end{equation}
When $\bb = 0$ our perturbation series agrees with the one known for $z_0^{\rm N}$.

\section{Large, or McMahon zeros}
To study large zeros we will use asymptotic expansions of Bessel functions as $z \to \infty$ \cite[\S9.2]{AbramowitzStegunHandbook2013}. If we denote 
$$\phi_{\nu}(z) = \sum_{k = 0}^{\infty}(-1)^{k}\frac{a_{2k}(\nu)}{z^{2k}}, \quad \psi_{\nu}(z) = \sum_{k = 0}^{\infty}(-1)^{k}\frac{a_{2k + 1}(\nu)}{z^{2k + 1}},$$
with
$$a_{0}(\nu) = 1, \quad a_{k}(\nu) = \frac{(4\nu^2 - 1)(4\nu^2 - 3)\ldots(4\nu^2 - (2k - 1)^{2})}{k! 8^{k}},$$
then
$$J_{\nu}(z) \sim \sqrt{\frac{2}{\pi z}} \left( \phi_{\nu}(z) \cos\omega - \psi_{\nu}(z) \sin\omega \right), \quad Y_{\nu}(z) \sim \sqrt{\frac{2}{\pi z}} \left( \phi_{\nu}(z) \sin\omega + \psi_{\nu}(z) \cos\omega \right),$$
where $\omega = z - \frac{\pi\nu}{2} - \frac{\pi}{4}$.
Similarly, functions
$$\tilde{\phi}_{\nu}(z) = \sum_{k = 0}^{\infty}(-1)^{k}\frac{b_{2k}(\nu)}{z^{2k}}, \quad \tilde{\psi}_{\nu}(z) = \sum_{k = 0}^{\infty}(-1)^{k}\frac{b_{2k + 1}(\nu)}{z^{2k + 1}},$$
where
\begin{align*}
b_{0}(\nu) &= 1, \quad b_1(\nu) = \frac{4\nu^2 + 3}{8}, \\
b_{k}(\nu) &= \frac{(4\nu^2 - 1)(4\nu^2 - 3)\ldots(4\nu^2 - (2k - 3)^{2}) \cdot (4\nu^2 + 4k^2 - 1)}{k! 8^{k}},
\end{align*}
produce
$$J'_{\nu}(z) \sim -\sqrt{\frac{2}{\pi z}} \left( \tilde{\phi}_{\nu}(z) \sin\omega + \tilde{\psi}_{\nu}(z) \cos\omega \right), \quad Y'_{\nu}(z) \sim \sqrt{\frac{2}{\pi z}} \left( \tilde{\phi}_{\nu}(z) \cos\omega - \tilde{\psi}_{\nu}(z) \sin\omega \right).$$
Denote $\mathcal{S} = \sin\left((\kk-1)z\right)$ and $\mathcal{C} = \cos\left((\kk-1)z\right)$. The cross-products can be expressed as
\begin{equation*}
\begin{split}
\frac{\sqrt{\kk}\pi z}{2} G^{1,1}_{\nu}(\kk, z) &\sim \Big[\tilde{\phi}_{\nu}(z)\tilde{\phi}_{\nu}(\kk z) + \tilde{\psi}_{\nu}(z)\tilde{\psi}_{\nu}(\kk z)\Big]\mathcal{S} + \Big[\tilde{\phi}_{\nu}(z)\tilde{\psi}_{\nu}(\kk z) - \tilde{\psi}_{\nu}(z)\tilde{\phi}_{\nu}(\kk z)\Big]\mathcal{C},\\
\frac{\sqrt{\kk}\pi z}{2} G^{0,1}_{\nu}(\kk, z) &\sim \Big[\phi_{\nu}(z)\tilde{\phi}_{\nu}(\kk z) + \psi_{\nu}(z)\tilde{\psi}_{\nu}(\kk z)\Big]\mathcal{C} - \Big[\phi_{\nu}(z)\tilde{\psi}_{\nu}(\kk z) - \psi_{\nu}(z)\tilde{\phi}_{\nu}(\kk z)\Big]\mathcal{S},\\
\frac{\sqrt{\kk}\pi z}{2} G^{1,0}_{\nu}(\kk, z) &\sim \Big[\tilde{\phi}_{\nu}(z)\psi_{\nu}(\kk z) - \tilde{\psi}_{\nu}(z)\phi_{\nu}(\kk z)\Big]\mathcal{S} - \Big[\tilde{\phi}_{\nu}(z)\phi_{\nu}(\kk z) + \tilde{\psi}_{\nu}(z)\psi_{\nu}(\kk z)\Big]\mathcal{C},\\
\frac{\sqrt{\kk}\pi z}{2} G^{0,0}_{\nu}(\kk, z) &\sim \Big[\phi_{\nu}(z)\psi_{\nu}(\kk z) - \psi_{\nu}(z)\phi_{\nu}(\kk z)\Big]\mathcal{C} + \Big[\phi_{\nu}(z)\phi_{\nu}(\kk z) + \psi_{\nu}(z)\psi_{\nu}(\kk z)\Big]\mathcal{S}.
\end{split}
\end{equation*}
Based on this we can evaluate
\begin{align*}
\frac{\sqrt{\kk}\pi z}{2} G^{1,1}_{\nu}(\kk, z) &= \bigg[ 1 + \bigg( - b_2 + \frac{b_1^2}{\kk} - \frac{b_2}{\kk^2}\bigg)\frac{1}{z^2} + \ldots \bigg]\mathcal{S} + \\
&+\bigg[ \bigg(-b_1 + \frac{b_1}{\kk}\bigg)\frac{1}{z} + \bigg(b_{3} - \frac{b_{1} b_{2}}{\kappa} + \frac{b_{1} b_{2}}{\kappa^{2}} - \frac{b_{3}}{\kappa^{3}}\bigg)\frac{1}{z^3}  + \ldots \bigg] \mathcal{C},\\
\frac{\sqrt{\kk}\pi z}{2} G^{0,1}_{\nu}(\kk, z) &= \bigg[ 1 + \bigg(- a_2 + \frac{a_1 b_1}{\kk} - \frac{b_2}{\kk^2}\bigg)\frac{1}{z^2} + \ldots \bigg]\mathcal{C} - \\
&-\bigg[ \bigg(-a_1 + \frac{b_1}{\kk}\bigg)\frac{1}{z} + \bigg(a_{3} - \frac{a_{2} b_{1}}{\kappa} + \frac{a_{1} b_{2}}{\kappa^{2}} - \frac{b_{3}}{\kappa^{3}}\bigg)\frac{1}{z^3}  + \ldots \bigg]\mathcal{S},
\end{align*}
\begin{align*}
\frac{\sqrt{\kk}\pi z}{2} G^{1,0}_{\nu}(\kk, z) &= \bigg[ \bigg(-b_1 + \frac{a_1}{\kk}\bigg)\frac{1}{z} +\bigg(b_{3} - \frac{a_{1} b_{2}}{\kappa} + \frac{a_{2} b_{1}}{\kappa^{2}} - \frac{a_{3}}{\kappa^{3}}\bigg)\frac{1}{z^3}  + \ldots \bigg]\mathcal{S} - \\
&-\bigg[ 1 + \bigg(- b_2 + \frac{a_1 b_1}{\kk} - \frac{a_2}{\kk^2}\bigg)\frac{1}{z^2} + \ldots \bigg]\mathcal{C},\\
\frac{\sqrt{\kk}\pi z}{2} G^{0,0}_{\nu}(\kk, z) &= \bigg[ \bigg(-a_1 + \frac{a_1}{\kk}\bigg)\frac{1}{z} + \bigg( a_{3} - \frac{a_{1} a_{2}}{\kappa} + \frac{a_{1} a_{2}}{\kappa^{2}} - \frac{a_{3}}{\kappa^{3}}\bigg)\frac{1}{z^3}  + \ldots \bigg]\mathcal{C} + \\
&+\bigg[ 1 + \bigg(- a_2 + \frac{a_1^2}{\kk} - \frac{a_2}{\kk^2}\bigg)\frac{1}{z^2} + \ldots \bigg]\mathcal{S}.
\end{align*}
Once we introduce new coefficients $$c_k(\bb,\nu) = b_k(\nu) + i\bb\nu a_{k-1}(\nu), \quad k\in\N,$$ the function $g_{\nu}$ becomes
\begin{multline*}
\frac{\sqrt{\kk}\pi z}{2} g_{\nu}(\bb,\kk,z) = \left[ 1 + \left( \frac{c_1^{2}}{\kk} - \left(1 + \frac{1}{\kk^2}\right)c_2  \right)\frac{1}{z^2} + \ldots \right] \mathcal{S} + \\
+ \left[ \frac{(1-\kk)c_1}{\kk}\frac{1}{z} + \left( \frac{\kk^3 - 1}{\kk^3}c_3 - \frac{\kk - 1}{\kk^2}c_2c_1 \right)\frac{1}{z^3}  + \ldots \right] \mathcal{C}.
\end{multline*}
Let
\begin{equation}\label{tangent}
	\tan\theta = \frac{\frac{(1-\kk)c_1}{\kk}\frac{1}{z} + \left( \frac{\kk^3 - 1}{\kk^3}c_3 - \frac{\kk - 1}{\kk^2}c_2c_1 \right)\frac{1}{z^3}  + \ldots}{1 + \left( \frac{c_1^{2}}{\kk} - \left(1 + \frac{1}{\kk^2}\right)c_2  \right)\frac{1}{z^2} + \ldots}
\end{equation}
then
\begin{equation*}
\frac{\sqrt{\kk}\pi z}{2} g_{\nu}(\bb,\kk,z) = \sin\left( (\kk - 1)z + \theta \right)
\end{equation*}
and the zeros are given by
\begin{equation*}
z = \frac{s\pi - \theta}{\kk - 1}, \quad s \in \Z.
\end{equation*}
Time to employ McMahon series \cite{McMahonRoots1894}. To this end, we expand \eqref{tangent} into power series
$$\tan\theta = \frac{(1-\kk)c_1}{\kk}\frac{1}{z} + \left[ \frac{\kk - 1}{\kk^2}c_1^3 + \frac{\kk^3 - 1}{\kk^3}(c_3 - c_2 c_1) \right]\frac{1}{z^3} + \ldots$$
followed by the expansion of $\theta$ itself
$$\theta = \frac{(1-\kk)c_1}{\kk}\frac{1}{z} + \frac{\kk^3 - 1}{\kk^3}\left( \frac{c_1^3}{3} + c_3 - c_2c_1 \right)\frac{1}{z^3} + \ldots$$
We arrive at a transcendental equation
\begin{equation*}
	z = \frac{s\pi}{\kk - 1} + \frac{c_1}{\kk}\frac{1}{z} - \frac{\kk^3 - 1}{\kk^3(\kk - 1)}\left( \frac{c_1^3}{3} + c_3 - c_2c_1 \right)\frac{1}{z^3} + \ldots,
\end{equation*}
which---if $s \neq 0$, $p = c_1 / \kk$, and $q = - \frac{\kk^3 - 1}{\kk^3(\kk - 1)}\left( \frac{c_1^3}{3} + c_3 - c_2c_1 \right)$---McMahon solves as
\begin{equation}\label{nonexceptional}
	z = \frac{s\pi}{\kk - 1} + p\left( \frac{s\pi}{\kk - 1} \right)^{-1} + (q - p^2)\left( \frac{s\pi}{\kk - 1} \right)^{-3} + \ldots.
\end{equation}
When $\bb = 0$, we have
$$p = \frac{4\nu^2 + 3}{8\kk}, \quad q = -\frac{\kk^3 - 1}{\kk^3(\kk - 1)}\left( \frac{b_1^3}{3} + b_3 - b_2b_1 \right) = \frac{\kk^3 - 1}{\kk^3(\kk - 1)}\frac{16\nu^4 + 184\nu^2 - 63}{384}$$
which coincides with the known expansion for the Neumann problem.

\section{Discussion}
We found asymptotic expressions \eqref{exceptional} and \eqref{nonexceptional} for positive-real-part-zeros of \eqref{g} as $\kk \to 1$
\begin{align*}
z_0(\bb) &= \nu\sqrt{\bb^2 + 1} \left[ 1 - \frac{\kk - 1}{2} - \left( \frac{5\bb^2 - 7}{\bb^2 + 1} + 4 i \bb \nu \right)\frac{(\kk - 1)^2}{24} + \ldots \right],\\
z_s(\bb) &= \frac{s\pi}{\kk - 1} + \frac{ \frac{4\nu^2 + 3}{8} + i\bb\nu }{\kk}\left(\frac{s\pi}{\kk - 1}\right)^{-1} + \ldots, \quad s \in \N
\end{align*}
These expressions hold for a fixed $\bb$, turn into $z_s^{\rm N}$ when $\bb = 0$, and can be seen as their branches. We have not shown that $z_s(\bb)$ exhaust all the zeros of \eqref{g} but it could, perhaps, be done in the same vein as was done by Cochran for $z_s^{\rm N}$ \cite{CochranRemarks1964}.

\begin{figure}[ht]
\centering
	\includegraphics[width=0.9\linewidth]{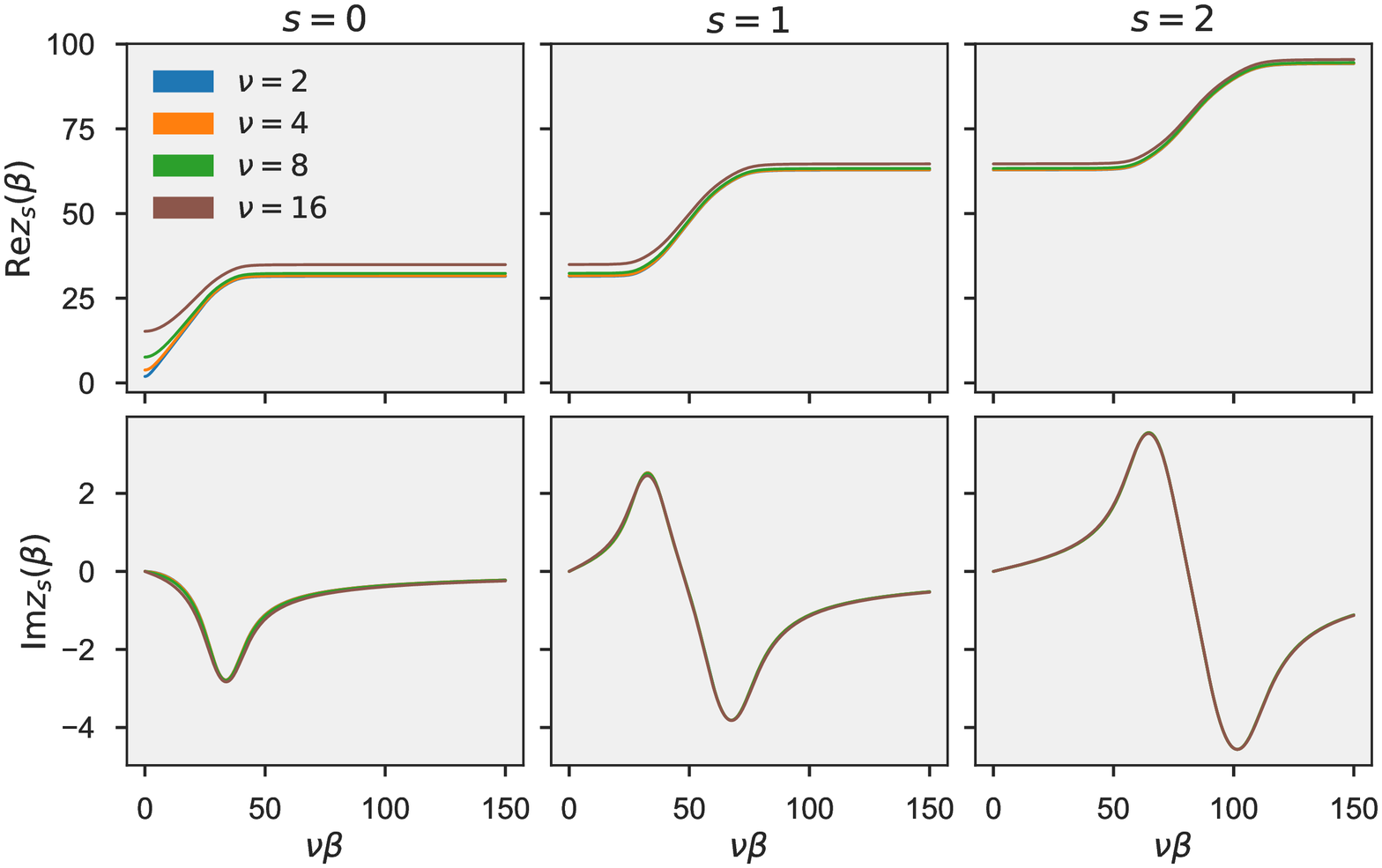}
    \caption{Numerically computed zeros $z_0(\bb)$, $z_1(\bb)$, and $z_2(\bb)$ for $\kk = 1.1$ and $\nu = 2,4,8,16$.}
    \label{fig:zeros}
\end{figure}

Now, what if we let $\bb$ grow? We used parameter continuation starting from $z_s^{\rm N}$ to evaluate $z_s(\bb)$ for large values of $\bb$. In fact we used $\nu\bb$ as the parameter rather than $\bb$. Numerical results are shown in Figure \ref{fig:zeros}.

One of the foremost properties to observe is that
$$z_s(\bb) \longrightarrow z_{s + 1}^{\rm D}, \quad \bb \to +\infty.$$
For instance, the "exceptional" zero $z_0^{\rm N}$ of the Neumann problem shapeshifts into the first McMahon zero $z_1^{\rm D}$ of the Dirichlet problem.

Next on the list is {\it how} they reach the limit. The real parts do so in finite time: when $s = 0$ it surges towards $\mathrm{Re}z_{1}^{\rm D}$ right from the start, while for $s > 0$ the real parts $\mathrm{Re}z_s(\bb)$ have a period of idleness before they undergo a phase transition with a plateau at $\mathrm{Re}z_{s+1}^{\rm D}$. As for the imaginary parts---unlike their monotone increasing real counterparts---they manage to have extrema: $\mathrm{Im}z_0(\bb)$ has a unique negative local minimum and when $s > 0$ there is a positive maximum followed by a negative minimum. In both cases extrema are followed by monotone $\bb^{-1}$-like convergence to zero. It is also remarkable that for every $s \in \N$ there is a $0<\bb<\infty$ such that $z_s(\bb)$ is real.

\begin{figure}[ht]
\begin{subfigure}[b]{0.5\textwidth}
\centering
	\includegraphics[width=0.9\textwidth]{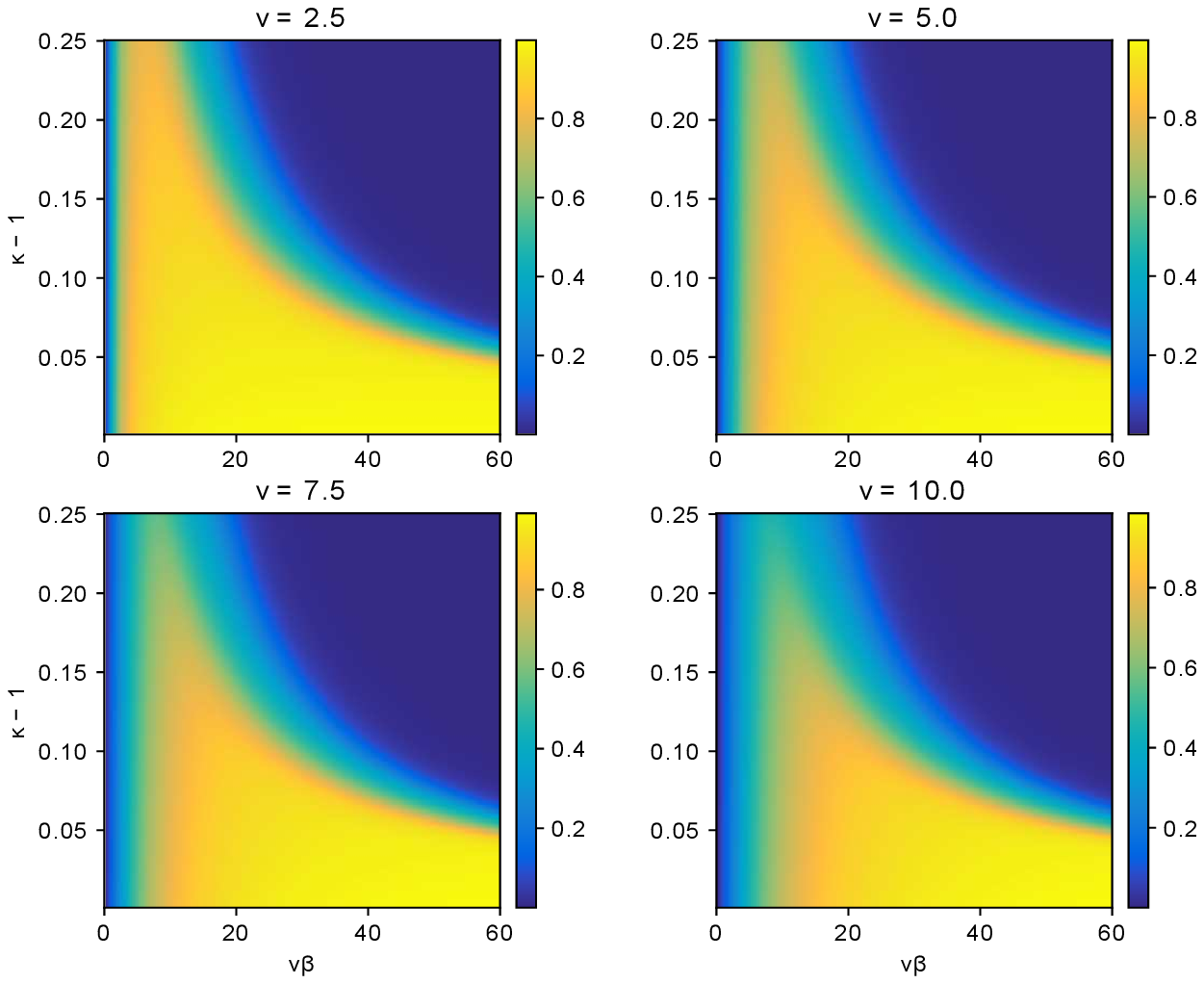}
    \caption{$| \mathrm{Re} \, \partial z_0 / \partial(\nu\beta) |$}
\end{subfigure}
\begin{subfigure}[b]{0.5\textwidth}
\centering
	\includegraphics[width=0.922\textwidth]{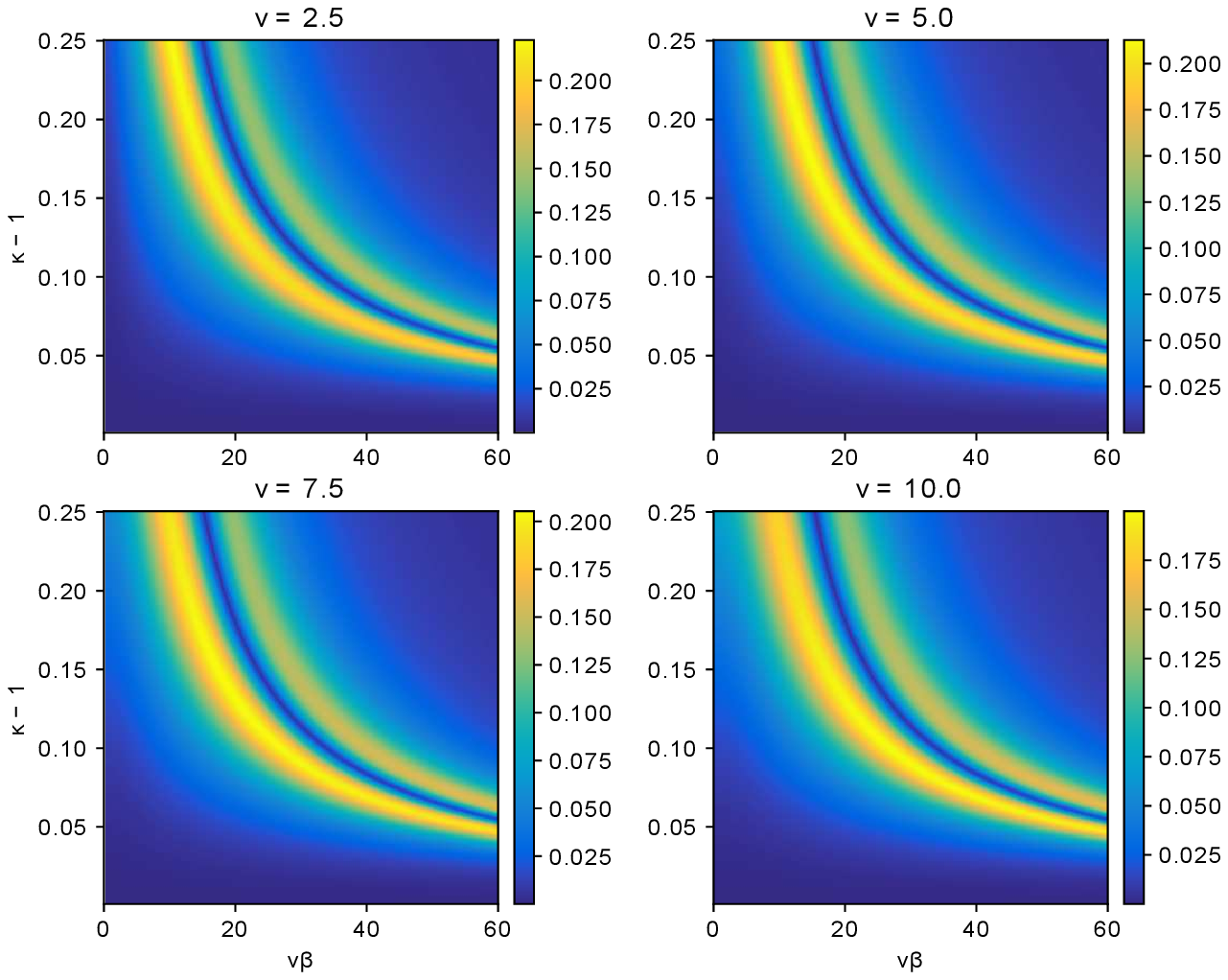}
    \caption{$| \mathrm{Im}\,\partial z_0 / \partial(\nu\beta) |$}
\end{subfigure}
    \caption{Derivative of $z_0(\bb)$ with respect to $\nu\beta$. The moment $\mathrm{Re}z_0$ reaches the plateau and the minimum of $\mathrm{Im}z_0$ correspond to zeros of the derivative and stand out as hyperbolas on the plot.}
    \label{fig:grads0}
\end{figure}

\begin{figure}[ht]
\begin{subfigure}[b]{0.5\textwidth}
\centering
	\includegraphics[width=0.9\textwidth]{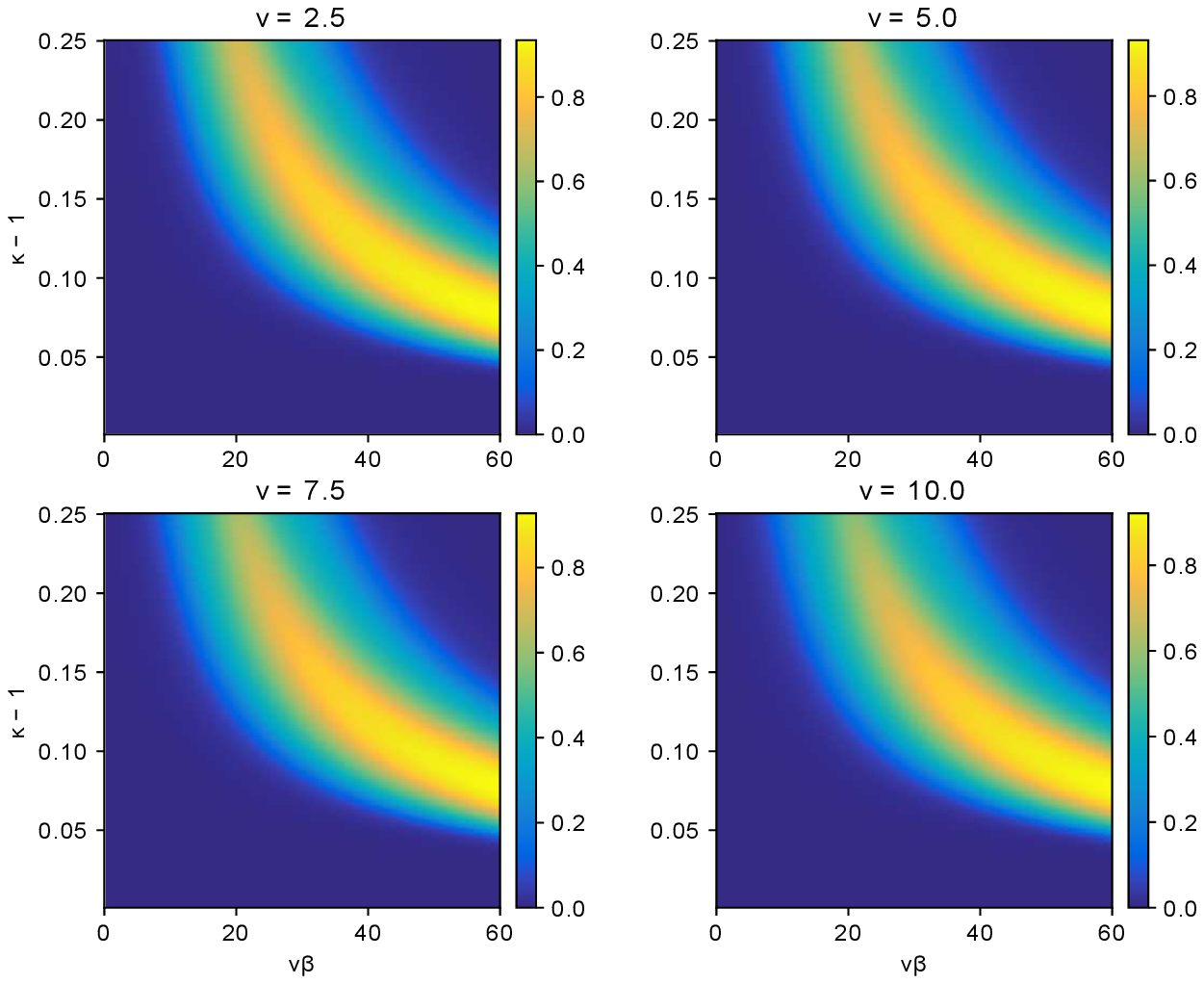}
    \caption{$| \mathrm{Re} \, \partial z_1 / \partial(\nu\beta) |$}
\end{subfigure}
\begin{subfigure}[b]{0.5\textwidth}
\centering
	\includegraphics[width=0.922\textwidth]{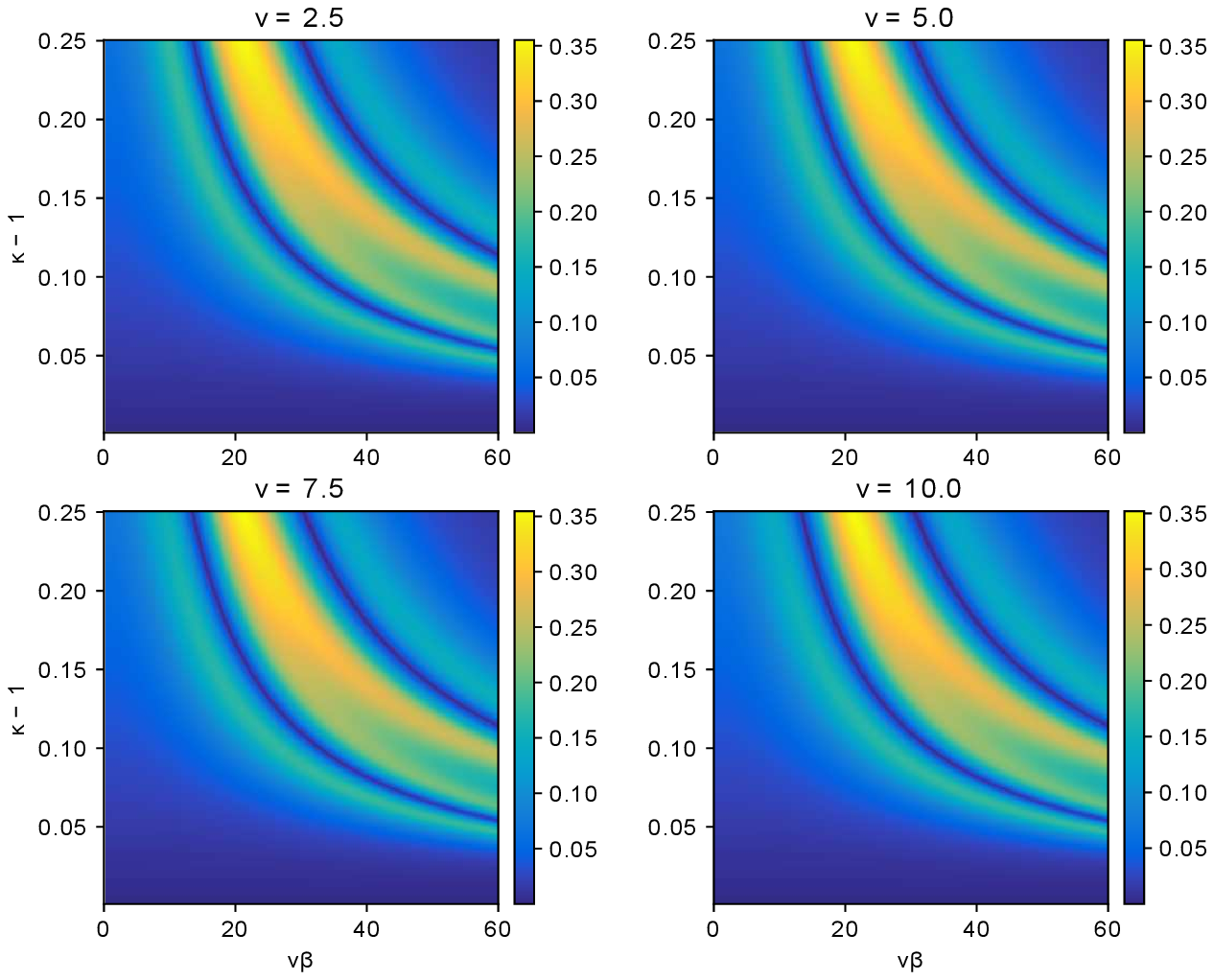}
    \caption{$| \mathrm{Im}\,\partial z_1 / \partial(\nu\beta) |$}
\end{subfigure}
    \caption{Derivative of $z_1(\bb)$ with respect to $\nu\beta$. The phase transition region of $\mathrm{Re}z_1$ and the extrema of $\mathrm{Im}z_1$ correspond to zeros of the derivative and stand out as hyperbolas on the plot.}
    \label{fig:grads1}
\end{figure}

To get some insight into the qualitative behaviour of $z_s(\bb)$, we plot $$\left| \mathrm{Re} \frac{\partial z_s}{\partial (\nu\beta)}  \right| \,\text{and}\, \left| \mathrm{Im} \frac{\partial z_s}{\partial (\nu\beta)}  \right|$$
against $\kappa - 1$ and $\nu\beta$ in Figures \ref{fig:grads0} and \ref{fig:grads1}. The salient points of the graph of $z_s(\bb)$---boundaries of the phase transition regions of $\mathrm{Re} z_s$ and extrema of $\mathrm{Im} z_s$---emerge as hyperbolas on the plots. In other words, the behaviour of $z_s(\bb)$ experiences qualitative changes when $\nu\beta(\kk - 1)$ reaches certain critical values.

The origins of the function $g_\nu$ shed some light onto the significance of $\nu\beta(\kk - 1)$. Recall that it comes from an oblique derivative boundary value problem for the Laplacian in an annulus with inner radius 1 and outer radius $\kappa$; $\beta$ is the tangent of the oblique angle. The order $\nu$ in this case is an integer and has the meaning of angular frequency. Let $1 \leqslant\rho \leqslant \kk$ and $0 \leqslant \theta <2\pi$ be polar radius and angle, respectively. Consider a logarithmic spiral 
$$\rho = e^{\frac{\theta - \theta_0}{\nu\bb}}$$
emanated from the inner circle at $\theta_0$. It reaches the outer circle at $$\theta = \theta_0 + \nu\bb \log{\kk} \approx \theta_0 + \nu\bb (\kk - 1)$$
and we see that for thin annuli the phase shift is approximately equal to $\nu\beta(\kk - 1)$.

As a side note, $\log\kk$ appears naturally when one studies imaginary $\nu$-zeros \cite{CochranRemarks1964,MARTINEKEtAlZeros1966}.
\section{Conclusion}
In this paper we looked at a combination of Bessel cross-products that comes from an oblique derivative boundary value problem in an annulus and derived asymptotic expressions for its zeros as the annulus collapses into a circle. We found that just like in the Neumann problem, there are two kinds of zeros: those that blow up and those that stay finite. To study the behaviour of the zeros as the oblique parameter $\bb \to +\infty$ (i.e. as the oblique derivative tends to the tangential one) we computed them numerically and observed that they experience qualitative changes---such as phase transitions or local extrema---when $\nu\bb(\kk-1)$ passes critical values.

\paragraph{Acknowledgements} The reported study was funded by RFBR according to the research project 18-31-00236.

\bibliographystyle{unsrt}
\bibliography{bib}
\end{document}